\documentclass[paper]{geophysics}

\usepackage{algorithm}
\usepackage{algorithmic}
\usepackage{url}
\usepackage{color}

\newcommand{\bfm}{\mathbf{m}}
\newcommand{\bfG}{\mathbf{G}}
\newcommand{\bfd}{\mathbf{d}}

\newcommand{\bfr}{\mathbf{r}}

\newcommand{\bfh}{\mathbf{h}}

\newcommand{\bfui}{\mathbf{u}_i}
\newcommand{\bfvi}{\mathbf{v}_i}

\newcommand{\WL}{\mathbf{W_{{L}_1}}}
\newcommand{\Wh}{\mathbf{W_{h}}}
\newcommand{\U}{\mathbf{U}}
\newcommand{\V}{\mathbf{V}}
\newcommand{\B}{\mathbf{B}}
\newcommand{\Si}{\mathbf{\Sigma}}
\newcommand{\D}{\mathbf{D}}
\newcommand{\W}{\mathbf{W}}
\newcommand{\bfQ}{\mathbf{Q}}

\newcommand{\bfdo}{\mathbf{d}_{\mathrm{obs}}}

\newcommand{\bfde}{\mathbf{d}_{\mathrm{exact}}}
\newcommand{\bfma}{\mathbf{m}_{\mathrm{apr}}}

\newcommand{\Wd}{\mathbf{W_{\bfd}}}

\newcommand{\Gtildetilde}{\tilde{\tilde{\bfG}}}

\newcommand{\Wz}{\mathbf{W_{\mathrm{z}}}}

\newcommand{\diag}{\mathrm{diag}}

\newcommand{\Rn}{\mathcal{R}^{n}}
\newcommand{\Rm}{\mathcal{R}^{m}}
\newcommand{\Rmn}{\mathcal{R}^{m \times n}}
\newcommand{\Rmq}{\mathcal{R}^{m \times q}}
\newcommand{\Rqq}{\mathcal{R}^{q \times q}}
\newcommand{\Rnq}{\mathcal{R}^{n \times q}}
\newcommand{\bfeta}{\mbox{\boldmath{$\eta$}}}

\begin{document}

\title{A fast algorithm for regularized focused  $3$-D inversion of  gravity data using the randomized SVD}

\renewcommand{\thefootnote}{\fnsymbol{footnote}} 


\address{
\footnotemark[1]Institute of Geophysics, \\
University of Tehran, \\
Tehran, Iran. Email: svatan@ut.ac.ir \\
\footnotemark[2]School of Mathematical and Statistical Sciences, \\
Arizona State University, \\
Tempe, AZ, USA. Email: renaut@asu.edu \\
\footnotemark[3]Institute of Geophysics, \\
University of Tehran, \\
Tehran, Iran. Email: ebrahimz@ut.ac.ir }

\author{Saeed Vatankhah\footnotemark[1], Rosemary Anne Renaut\footnotemark[2] and Vahid Ebrahimzadeh Ardestani\footnotemark[3] }

\lefthead{Vatankhah et al.}
\righthead{Inversion of gravity data using the RSVD}

\maketitle

\begin{abstract}
A fast algorithm for solving the under-determined $3$-D linear gravity inverse problem based on the randomized singular value decomposition (RSVD) is developed. The algorithm combines an iteratively reweighted approach for $L_1$-norm regularization with the RSVD methodology in which the large scale linear system at each iteration is replaced with a much smaller linear system.  Although the optimal choice for the low rank approximation of the system matrix with $m$ rows is $q=m$, acceptable results are achievable with $q\ll m$.  In contrast to the use of the LSQR algorithm for the solution of the linear systems at each iteration, the singular values generated using the RSVD yield a good approximation of the dominant singular values of the large scale system matrix. The regularization parameter found for the small system at each iteration is thus dependent on the dominant singular values of the large scale system matrix and appropriately regularizes the dominant singular space of the large scale problem. The results achieved are comparable with those obtained using the LSQR algorithm for solving each linear system, but are obtained at reduced computational cost.  The method has been tested on synthetic models along with the real gravity data from the Morro do Engenho complex from central Brazil.
\end{abstract}

\section{Introduction}

It is well-known that the linear gravity inverse problem is ill-posed, and that effective regularization methods should be used to obtain reasonable solutions, \citep{LiOl:98,PoZh:99,BoCh:2001,SiBa:2006}.  In geophysical inverse modeling it is often assumed that the sources of interest are localized and separated by distinct interfaces. Thus, the inversion methodology should be able to provide sharp and focused images of the subsurface. Many different approaches have been used, including  the compactness constraint \citep{LaKu:83}, minimum gradient support \citep{PoZh:99,Zhdanov:2002}, total variation regularization \citep{BCO:2002},  applying the Cauchy norm on the model parameters \citep{Pi:09} and using the $L_1$-norm stabilizer \citep{Far:2008,Loke:2003,VRA:2017}. In all these cases, the process is iterative and  the model-space iteratively reweighted least squares (IRLS) algorithm may be used. Here, we suppose a focused image of the subsurface is preferred and  adopt  the $L_1$ inversion methodology presented in Vatankhah et. al \shortcite{VRA:2017} for determining the solution of the under-determined inversion problem with $m$ data measurements for recovery of a volume with $n$ cells, $m \ll n$.  

For $G \in \mathcal{R}^{m \times n}$ with both $m$ and $n$ relatively small, a physically acceptable numerical solution is obtained using the singular value decomposition (SVD), or the generalized singular value decomposition (GSVD), as appropriate. For large-scale inverse problems it is no longer feasible, whether with respect to memory or computational time, to rely on a direct solver \citep{OlLi:94,LiOl:03}. Rather, nowadays, the LSQR algorithm based on the Golub-Kahan bidiagonalization (GKB) is frequently used   \citep{PaSa:1982a, PaSa:1982b,KiOl:2001,ChNaOl:2008,RVA:2017,VMN:2015,VRA:2017}. For the LSQR algorithm using $t \ll m$ steps of the GKB process a Krylov subspace with dimension $t$ is generated and  the solution is obtained on this subspace at negligible computational cost  using the SVD of the subspace system matrix.  On the other hand, randomized algorithms can be used to efficiently and directly approximate the SVD of $G$ \citep{Halko:2011,XiZo:2013,VMN:2015} yielding a rank $q$ approximation of $G$  in which $q\ll m$. Here, we employ the randomized SVD algorithm with Gaussian sampling  which allows us to compute the rank $q$ SVD approximation of $G$ efficiently.

For both LSQR and RSVD algorithms a suitable value for $t$ and $q$, respectively, must be determined. Generally both $t$ and $q$ should be as small as possible in order that the inversion methodology is fast. For the LSQR algorithm it is known that it is important that the choice for $t$  provides an approximate system matrix that accurately capture the dominant spectral properties of the original system matrix, and that  an optimal regularization parameter can be found providing effective regularization of the large-scale problem,  \citep{KiOl:2001,ChNaOl:2008,RVA:2017}. Here we demonstrate an approach for selecting $q$ dependent on $m$ which simultaneously yields effective system matrix approximation  and appropriate regularization parameter estimation. We contrast the RSVD technique  with the inversion methodology based on the LSQR algorithm, Vatankhah et al. \shortcite{VRA:2017} and demonstrate  that  solutions of comparable quality are obtained at reduced computational cost. Furthermore, the method of unbiased predictive risk estimation (UPRE) for finding the regularization parameter has to be modified to use a truncated spectrum when applied in the context of the LSQR algorithm \citep{VRA:2017} but can be used directly with the RSVD methodology. Consequently,   we present a fast methodology for inversion of gravity data with $L_1$ regularization using a new RSVD methodology.

\section{Inversion methodology}
 
We will consider the under-determined linear system arising in inversion of gravity data, see Li \& Oldenburg \shortcite{LiOl:98} and Boulanger \& Chouteau \shortcite{BoCh:2001},
\begin{equation}\label{d=gm}
\bfdo= \bfG \bfm.
\end{equation}
$\bfG \in \Rmn$  is the forward modeling operator, and vectors $\bfdo \in \Rm$ and $\bfm \in \Rn$ contain noisy measurement data and unknown model parameters, the densities of cells, respectively. The goal is to find a geologically acceptable model which satisfies the observed data at the noise level. The problem is ill-posed and  regularization is required to achieve a meaningful solution. We use the $L_1$-norm regularization methodology presented in Vatankhah et. al \shortcite{VRA:2017} in which the solution of~(\ref{d=gm}) is obtained from the minimization of the following non-linear objective function, 
\begin{equation}\label{globalfunction1}
P^{\alpha}(\bfm)=\| \Wd(\bfG\bfm-\bfdo)  \|_2^2 + \alpha^2 \| \W(\bfm-\bfma)\|_2^2.
\end{equation}
The matrix $\W$ is the product of three diagonal matrices, a depth weighting matrix $\Wz$, \cite{LiOl:98}, a matrix $\WL={((\bfm-\bfma)^2+\epsilon^2)^{-1/4}}$ arising from approximation of the $L_1$-norm stabilizer with a $L_2$-norm term, and a hard constraint matrix $\Wh$. The vector $\bfma$ is either set to zero or  is an initial model selected based on prior available information, \citep{LiOl:96}.  If the densities of some cells are known they are used in $\bfma$ and  a large value for their corresponding entries on the diagonal of $\Wh$ is selected. Otherwise the entries in $\Wh$ are $1$. The matrix $\Wd^{-1}=\diag (\bfeta)$  is a data weighting matrix in which the component $\eta_i$ of $\bfeta$ is the standard deviation of the noise in the $i$th datum. The regularization parameter $\alpha$ balances two terms in objective function~(\ref{globalfunction1}), and its determination is an important step in any regularization method. Note that  the minimum support constraint can be used in~(\ref{globalfunction1}) by simply replacing $\WL$ with $\mathbf{W_{\mathrm{MS}}}={((\bfm-\bfma)^2+\epsilon^2)^{-1/2}}$.

The matrix $\W$ is diagonal and  the objective function in (\ref{globalfunction1}) is easily transformed to the standard Tikhonov form, see Vatankhah et al. \shortcite{VAR:2015,VRA:2017}, as
\begin{equation}\label{globalfunction2}
P^{\alpha}(\bfh)=\| \Gtildetilde \bfh- \tilde{\bfr} \|_2^2 + \alpha^2 \|\bfh \|_2^2,
\end{equation}
with system matrix $\Gtildetilde=\tilde{\bfG}\W^{-1}=\Wd \bfG \W^{-1}$, right hand side residual vector $\tilde{\bfr}=\Wd (\bfdo-\bfG\bfma)$ and unknown model increment $\bfh=\W(\bfm-\bfma)$. The minimization of~(\ref{globalfunction2}) yields the solution
\begin{equation}\label{hsolution}
\bfh(\alpha)=((\Gtildetilde)^T\Gtildetilde+ \alpha^2 I_n)^{-1} (\Gtildetilde)^T \tilde{\bfr},
\end{equation}
and the model update 
\begin{equation}\label{modelupdate}
\bfm(\alpha)=\bfma+\W^{-1}\bfh (\alpha).
\end{equation}

The non-linearity in (\ref{globalfunction1}) arises because $\WL$, and hence $\W$, depends on the model parameters. We use the  model-space iteratively reweighted least squares (IRLS) algorithm to find the solution, as detailed  in Algorithm~\ref{IterativeL1RSVD}. The iteration is terminated when the  solution satisfies the noise level, $\chi_{\mathrm{Computed}}^2 =  \|\Wd {(\bfdo -\bfG\bfm)}\|_{2}^2 \leq m+\sqrt{2m} $, or a predefined maximum number of iterations, $K_{\mathrm{max}}$, is reached. Furthermore, the positivity constraint [$\rho_{\mathrm{min}}, \rho_{\mathrm{max}}$] is used to recover a reliable subsurface model. If at any iteration a density value falls outside these predefined lower and upper density bounds, the value is projected back to the nearest bound value \citep{BoCh:2001}.

For small $m$ and $n$, the solution $\bfh(\alpha)$ can be calculated cheaply using the  SVD: $U\Sigma V^T=\Gtildetilde$, see Algorithm~1 in Vatankhah et al. \shortcite{VAR:2015,VRA:2017}. Furthermore,  the availability of the SVD makes it possible to put the methods for regularization parameter estimation into  convenient and easy to use forms \citep{ChPa:2015,XiZo:2013}. It may be infeasible or expensive, however, to compute the SVD for large under-determined systems $m \ll n$. Even though we would only need the first $m$ columns of $V$, namely only the thin SVD, the cost is approximately $6nm^2+20m^3$, Golub and Van Loan \shortcite{GoVL:2013} page $493$.
In this case, a randomized algorithm can be used to approximate the SVD of $\Gtildetilde$. The method  uses random sampling to construct a low-dimensional subspace that captures most of the spectral properties  of the matrix, and then restricts the matrix to this subspace \citep{Halko:2011}. A standard factorization such as the SVD or eigen-decomposition can be applied  for the reduced matrix. Here, we develop a randomized SVD algorithm with Gaussian sampling for under-determined problems, see Algorithm~\ref{RSVDAlgorithm}, based on a combination of the methodologies presented in Voronin et al. \shortcite{VMN:2015} and Xiang and Zou \shortcite{XiZo:2013}. 
\begin{algorithm}
\caption{RSVD algorithm. Given $\Gtildetilde \in \Rmn$ ($m < n$), a target matrix rank $q $ and 
a small constant oversampling parameter $p$ satisfying $q+p=l  \ll m$, compute a low-rank approximation of $\Gtildetilde$: $\Gtildetilde_q = \U_q \Si_q \V_q^T$ with $\U_q \in \Rmq$, $\Si_q \in \Rqq$ and $\V_q \in \Rnq$.}\label{RSVDAlgorithm}
\begin{algorithmic}[1]
\STATE Generate a Gaussian random matrix $\mathbf{\Omega} \in \mathcal{R}^{l \times m} $.
\STATE Compute matrix $\mathbf{Y}=\mathbf{\Omega} \Gtildetilde \in \mathcal{R}^{l \times n}$.  
\STATE Compute orthonormal matrix $\bfQ \in \mathcal{R}^{n \times l}$ via QR factorization $\mathbf{Y}^T=\mathbf{Q}\mathbf{R}$.(Note that $Q$ is stored in factored form and not accumulated).
\STATE Form the matrix $\B=\Gtildetilde\bfQ \in \mathcal{R}^{m \times l}$ using factored form of $Q$.
\STATE Compute the matrix $\B^T\B \in \mathcal{R}^{l \times l}$.
\STATE Compute the eigen-decomposition of $\B^T\B$; $[\tilde{\V}_l, \D_l]=\mathrm{eig}(\B^T\B)$. 
\STATE Compute $\V_q=\bfQ \tilde{\V}_l(:,1:q)$;  $\Si_q= \sqrt{\D_l}(1:q,1:q)$; $ \mathrm{and}$ $ \U_q=\B \tilde{\V}_q(:,1:q) \Si_{q}^{-1}. $
\STATE Note $\tilde{\tilde{\bfG}}_q = \U_q\Si_q \V_q^T$
\end{algorithmic}
\end{algorithm}
 
In Algorithm~\ref{RSVDAlgorithm}, $p$ is a small oversampling parameter and provides a flexibility that is crucial for the effectiveness of the randomized SVD methods. We use the fixed value $p=10$ \citep{Halko:2011}. For clarification we explain the steps in Algorithm~\ref{RSVDAlgorithm}. Step~$2$ is used to extract the {range} of matrix $ (\Gtildetilde)^T$, i.e. $\mathrm{range}~\mathbf{Y}^T \approx \mathrm{range}~(\Gtildetilde)^T$. In step $3$, an orthogonal matrix $\bfQ$ is formed to represent the range of $(\Gtildetilde)^T$, and gives an approximation of the right singular vectors of $\Gtildetilde$. The original matrix is projected onto a lower dimensional matrix $\B$ in step~$4$, in which the matrix $\B$ provides the information on the range of $\Gtildetilde$, or the range of the left singular vectors. Step~$6$ provides the eigen-decomposition of matrix $\B^T\B$. These eigenvalues and eigenvectors are then used, step~$7$, to compute the singular values and the singular vectors of the matrix $\B$, yielding the rank $q$ SVD approximation to $\Gtildetilde$. See Appendix~\ref{eigtosv} for more details. For small $q$, this methodology is very effective because the reduced matrix $\B^T\B$ can be easily constructed and its eigen-decomposition rapidly computed. 
In general for an ill-conditioned matrix forming $A^TA$ will square the condition number. We note that $\Gtildetilde$ is mildly ill-conditioned and so $\Gtildetilde^T\Gtildetilde$ has condition $(\sigma_1/\sigma_r)^2$ where $\sigma_r$ is the smallest non zero singular value of $\Gtildetilde$. When $q$ is not too large as compared to $m$, $q\ll r$  the condition of $B^TB$ is approximately $(\sigma_1/\sigma_q)^2$, which is acceptable for the gravity inversion problem. Generally when the original matrix has more severe ill-conditioning we would not recommend the approach of forming $B^TB $ to find the singular decomposition for $B$, rather one would need to compute the SVD of $B$ in Algorithm~\ref{RSVDAlgorithm}, see Xiang and Zou \shortcite{XiZo:2013}.

We present the cost of each step of the Algorithm~\ref{RSVDAlgorithm} in Table~\ref{costtab}, in which we use the notation that a dot product of length $n$ costs  $2n$ flops. When $l \ll m\ll n$, the dominant cost is $6lmn$ and occurs for steps~$2$ and $4$, noting that in step~3 we assume the factorization is calculated without accumulating the matrix $Q$, see for example Golub and Van Loan \shortcite{GoVL:2013}. Thus Algorithm~\ref{RSVDAlgorithm} presents a great advantage in efficiency for large problems. Note, the parameter $q$ should be selected such that the dominant spectral properties of the original matrix are captured so that the solution obtained from the RSVD is close to the solution obtained using all components of the  full SVD. We will discuss how to choose $q$ in the next section. The iterative inversion methodology using the RSVD is given in Algorithm~\ref{IterativeL1RSVD}.

\begin{table}
\begin{center}
\caption{For each step (row 1) the level 3 (cubic) costs (row 2) of Algorithm~\ref{RSVDAlgorithm}.}\label{costtab}
\begin{tabular}{c  c c c c c c c}
\hline
$2$&$3$&$4$&$5$&$6$&$7$  \\ \hline
$2lmn$ & $2l^2(n-l/3)$  & $4lmn$& $2l^2m$ & $O(l^3)$ & $lq(2l+3m)$ \\
\hline
\end{tabular}
\end{center}
\end{table}

The other well-known algorithm for solving large inverse problems is the GKB based LSQR algorithm  \citep{PaSa:1982a,PaSa:1982b,KiOl:2001,ChNaOl:2008,RVA:2017}. In this case the original problem is projected onto a Krylov subspace using $t$ steps of the GKB process  and the solution can be obtained on this subspace cheaply using the SVD of the projected matrix. The details of the application of the method for the gravity inverse problem are presented in Vatankhah et al. \shortcite{VRA:2017}, Algorithm~$2$, and are thus not repeated here. The role of parameter $t$ in the LSQR here is the equivalent to that of  $q$ for the RSVD algorithm. While it is important that  $t \ll m$ in order to yield an efficient  and fast algorithm for large problems, simultaneously, $t$ should be selected large enough so that the dominant spectral properties of the original problem are accurately captured. In Vatankhah et al. \shortcite{VRA:2017}, it was suggested that  $t \geq m/20$ is suitable for the gravity problem. In the following we use both algorithms on  synthetic examples, and show that the presented algorithm based on the RSVD is significantly faster than that using the LSQR. 

\begin{algorithm}
\caption{Iterative $L_1$ Inversion using RSVD}\label{IterativeL1RSVD}
\begin{algorithmic}[1]
\REQUIRE $\bfdo$, $\bfma$, $\bfG$, $\Wd$, $\Wh$, $\epsilon > 0$, $\rho_{\mathrm{min}}$, $\rho_{\mathrm{max}}$, $K_{\mathrm{max}}$
\STATE Calculate $\Wz$, $\W^{(1)}=\Wz \Wh$, $\tilde{\bfG}= \Wd \bfG$
\STATE Initialize $\bfm^{(0)}=\bfma$, $(\WL)^{(1)} =I$, $k=0$
\STATE Calculate $\tilde{\bfr}^{(1)}=\Wd(\bfdo-\bfG\bfm^{(0)})$, $(\Gtildetilde)^{(1)}=\tilde{\bfG}(\W^{(1)})^{-1} $ 
\WHILE {Not converged, noise level not satisfied, and $k<K_{\mathrm{max}}$} 
\STATE {$k=k+1$}
\STATE {\label{start}Find the RSVD using Algorithm~\ref{RSVDAlgorithm}: $\Gtildetilde \approx \U_q\Si_q \V_q^T$} 
\STATE {\label{stepalpha}Use regularization parameter estimation to find $\alpha^{(k)}$}
\STATE {\label{stop}Set $\bfh^{(k)}= \sum_{i=1}^{q} \frac{\sigma_i^2}{\sigma_i^2+(\alpha^{(k)})^2} \frac{\bfui ^T\tilde{\bfr}^{(k)}}{\sigma_i} \bfvi$}
\STATE {Set $\bfm^{(k)}=\bfm^{(k-1)}+ (\W^{(k)})^{-1}\bfh^{(k)}$}
\STATE {Impose constraint conditions on $\bfm^{(k)}$ to force $\rho_{\mathrm{min}}\le \bfm^{(k)} \le \rho_{\mathrm{max}}$}
\STATE {Test convergence and exit loop if converged}
\STATE {Calculate the residual $\tilde{\bfr}^{(k+1)}=\Wd(\bfdo-G\bfm^{(k)})$}
\STATE {\label{Wupdate}Set $(\WL)^{(k+1)} =\diag\left( \left((\bfm^{(k)}-\bfm^{(k-1)})^2+\epsilon^2 \right)^{-1/4}\right)$, and  $\W^{(k+1)}=(\WL)^{(k+1)}\W^{(1)}$}
\STATE {Calculate $(\Gtildetilde)^{(k+1)}=\tilde{\bfG}(\W^{(k+1)})^{-1} $}
\ENDWHILE
\ENSURE Solution $\rho=\bfm^{(k)}$. $K=k$.
\end{algorithmic}
\end{algorithm}
 
\subsection{Regularization parameter estimation}

The estimation of a suitable regularization parameter $\alpha$ is an important and critical step in Algorithm~\ref{IterativeL1RSVD}. Among the many well-known approaches we use the UPRE, for which the derivation for (\ref{globalfunction2}) is given, for example, in Vogel \shortcite{Vogel:2002}, and is not repeated here. Using the RSVD approximation for $\Gtildetilde$, the UPRE function to be minimized is given by
\begin{equation}\label{UPRERSVD}
U(\alpha)=\sum_{i=1}^{q}  \left( \frac{1}{\sigma_i^2 \alpha^{-2} + 1} \right)^2 \left(\bfui^T\tilde{\bfr} \right)^2 + 2 \left( \sum_{i=1}^{q} \frac{\sigma_i^2}{\sigma_i^2+\alpha^2}\right) - q.
\end{equation} 
Typically $\alpha_{opt}$ is found by evaluating (\ref{UPRERSVD}) on a range of $\alpha$, between minimum and maximum $\sigma_i$, and then the value which minimizes the function is selected as $\alpha_{opt}$. In step~\ref{stepalpha} of Algorithm~\ref{IterativeL1RSVD} we use (\ref{UPRERSVD}) to estimate the regularization parameter.

The UPRE method for the LSQR algorithm was developed by Renaut et al. \shortcite{RVA:2017} and  used by Vatankhah et al. \shortcite{VRA:2017} for the gravity inverse problem. It was demonstrated that the direct application of the UPRE does not provide an $\alpha$ which is optimal for the Krylov subspace of size $t$ after $t$ steps of the LSQR algorithm.  It is known that the spectrum of the system matrix for the Krylov subspace solution (the projected matrix)  inherits the ill-conditioning of the original large-scale problem, \cite{PaSa:1982a, PaSa:1982b}. Further, it is shown in \cite{VRA:2017} that the model matrix for the gravity inverse problem is mildly ill-conditioned. While  the dominant spectral values of the projected matrix give a good approximation to the dominant spectrum of $\Gtildetilde$, due to inheriting the ill-conditioning the projected matrix also possesses very small spectral values. It is these small spectral values which lead to the estimation of a regularization parameter that is underestimated in relation to the original problem if using (\ref{UPRERSVD}) directly. Vatankhah et al. \shortcite{VRA:2017}  demonstrated that if the UPRE function is calculated for a truncated spectrum of the projected matrix, denoted as TUPRE, hence ignoring the small singular values, good estimates for $\alpha_{opt}$ are obtained. The iterative $L_1$ algorithm using the GKB and the TUPRE method is presented in Algorithm~$2$ in Vatankhah et al. \shortcite{VRA:2017}. Here we will show that the singular values estimated from the RSVD do not inherit the ill-conditioning of the $\Gtildetilde$  and  there is no need to truncate  the spectrum in applying (\ref{UPRERSVD})

\section{Synthetic examples}

Two synthetic examples are used here to analyse the presented methodology in Algorithm~\ref{IterativeL1RSVD}. The goal is at first to estimate the accuracy of the RSVD method as compared with using the full SVD (FSVD), and then to consider the relative computational costs. The first example is a small model consisting of two cubes. The model is selected small so that it is possible to contrast the accuracy of the solutions using both the FSVD and the RSVD. The second model, which is larger and consists of multiple bodies, is used to study the computational costs as compared to the use of the LSQR algorithm with the TUPRE, as presented in Vatankhah et al. \shortcite{VRA:2017}. The following tests are preformed on a desktop computer with Intel Core i$7$-$4790$ $3.6$G and $16$GB RAM.

\subsection{Model consisting of two cubes}

The first model consists of two similar cubes with density contrast $1$~g~cm$^{-3}$ embedded in a homogeneous background, Figure~\ref{fig:figure1a}. The cubes have dimension $300$~m $\times$ $300$~m $\times$ $200$~m and start at depth $50$~m. The gravity data, $\bfde$, are generated at the $30 \times 20 = 600$ points on the surface with grid spacing $50$~m. Gaussian noise with standard deviation of $ (0.02~(\bfde)_i + 0.002~\| \bfde \|)$ was added to each exact data point, yielding a noisy data set $\bfdo$, Figure~\ref{fig:figure1b}. The subsurface volume  is discretized into $6000$ similar cubes with cell size  of $50$~m in each dimension. The resulting kernel matrix  $\bfG$ of size  $600 \times 6000$  can be used for both the RSVD and FSVD inversion methodologies, and it is feasible to calculate the spectrum for the full problem, and thus examine the approximation of the spectrum of the FSVD by that of the RSVD . For the inversion methodology presented in Algorithm~\ref{IterativeL1RSVD} we use $\bfma=\mathbf{0}$, $\rho_{\mathrm{min}}=0$~g~cm$^{-3}$, $\rho_{\mathrm{max}}=1$~g~cm$^{-3}$ and $K_{\mathrm{max}}=50$. The number of iterations $K$, the final regularization parameter $\alpha^{(K)}$, the relative error of the reconstructed model
\begin{equation}\label{RE}
RE^{(K)}=\frac{\|\bfm_{\mathrm{exact}}-\bfm ^{(K)} \|_2}{\|\bfm_{\mathrm{exact}} \|_2},
\end{equation}
and the computational costs are reported in Table~\ref{tab1}. We should note here, as suggested by Farquharson  and  Oldenburg \shortcite{FaOl:2004}, that it is important to use $\alpha$ large at the first iteration and we thus use equation ($19$) in Vatankhah et al. \shortcite{VRA:2017} for $\alpha^{(1)}$ but the UPRE for subsequent $\alpha^{(k)}$. The results of the inversion using the FSVD, Algorithm~1 in Vatankhah et al. \shortcite{VAR:2015,VRA:2017},  and the LSQR with the TUPRE, Algorithm~2 in Vatankhah et al. \shortcite{VRA:2017}, are also presented. From Table~\ref{tab1} we can see that for $q=m$, the RSVD and FSVD  lead to the same models. This indicates that the RSVD Algorithm~\ref{RSVDAlgorithm} is consistent. For $q=50$ the inversion algorithm terminates at $K_{\mathrm{max}}$, i.e. the noise level is not satisfied,  and the error of the reconstructed model is large. With increasing $q$ the solution improves and achieves the accuracy of the FSVD algorithm. The inversion methodology using the LSQR with the TUPRE  algorithm has a different behavior. While acceptable solutions are obtained even with small $t$, and the accuracy improves with increasing $t$, the computational time increases dramatically. We note, for this example, the goal is not to compare the CPU time of the methods, but to verify the accuracy of the RSVD algorithm. Indeed the FSVD for small problems will always be faster than the use of the RSVD or LSQR algorithms, which both involve a first step of finding a subspace for the solution and then the generation of a low rank SVD. To illustrate the results using Algorithm~\ref{IterativeL1RSVD}, we present the reconstructed model for case $q=100$ in Figure~\ref{fig:figure2a}. As evident from the data given in Table~\ref{tab1} the results are acceptable and the error decreases with increasing $q$. The regularization parameter and relative error at each iteration are presented in Figures~\ref{fig:figure2b}-\ref{fig:figure2c}, respectively, and the UPRE functional at iteration $5$ in Figure~\ref{fig:figure2d}. The results demonstrate that, generally, Algorithm~\ref{IterativeL1RSVD} is able to reconstruct a sharp and focused image of the subsurface. 

\multiplot{2}{figure1a,figure1b}{width=0.4\textwidth}
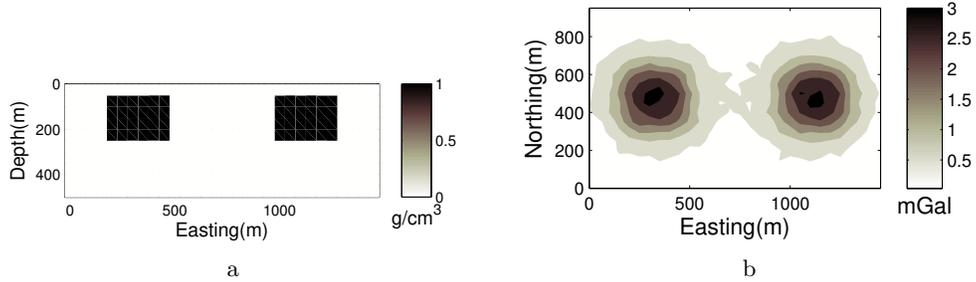
{(a) A model consisting of two cubes with density contrast $1$~g~cm$^{-3}$ embedded in a homogeneous background, cross-section at $y=475$~m; (b) The noisy gravity data of the model generated at the surface.}

\begin{table}
\begin{center}
\caption{The inversion results obtain from inversion methodologies: FSVD, LSQR and RSVD, using UPRE for FSVD and RSVD but TUPRE for LSQR. Different values of parameters $q$ and $t$ are used for LSQR and RSVD, respectively.}\label{tab1}
\begin{tabular}{c  c  c  c  c c}
\hline
Method&  $RE^{(K)}$&   $\alpha^{(1)}$ & $\alpha^{(K)}$& $K$ & Time (s)  \\ \hline
FSVD   & 0.3276& 55433& 53.52& 8& 15.8 \\ \hline
LSQR   &  & &  & &\\
$t=50$ & 0.4045& 16751& 31.29& 5& 10.6 \\ 
$t=100$ & 0.3240& 24247& 44.68& 7& 15.6\\ 
$t=150$ & 0.3142& 29944& 41.65& 7& 18.6\\ 
$t=200$ & 0.3097& 34681& 40.68& 7& 21.7\\ 
$t=600$ & 0.3262& 55433& 53.06& 8& 77.3\\  \hline
RSVD &  & &  & &\\
$q=50$ &0.4475& 18647& 42.78& 50& 83.3\\ 
$q=100$ & 0.3742& 26382& 56.40& 10& 17.6\\ 
$q=150$ & 0.3467& 31771& 47.68& 9& 16.4\\ 
$q=200$ & 0.3425& 36057& 50.50& 9& 16.2\\ 
$q=600$ & 0.3276& 55433& 53.52& 8& 18.3\\  \hline
\end{tabular}
\end{center}
\end{table}

\multiplot{4}{figure2a,figure2b,figure2c,figure2d}{width=0.4\textwidth}
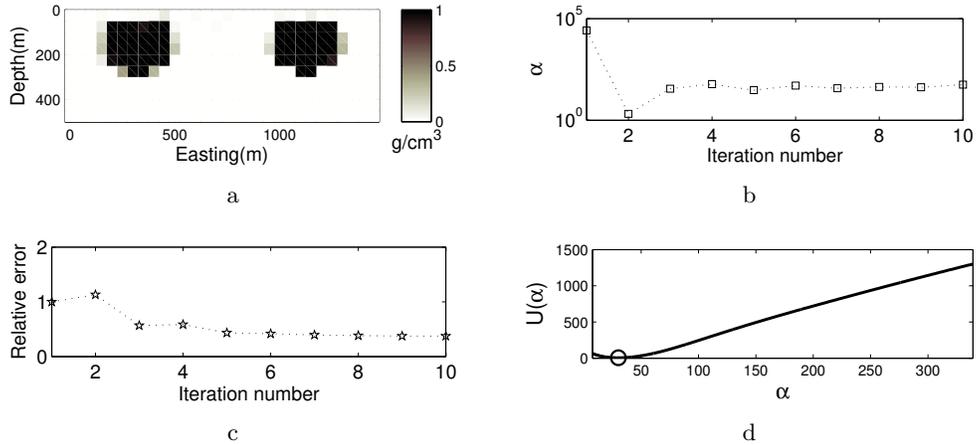
{The results of the inversion using Algorithm~\ref{IterativeL1RSVD} when $q=100$ is selected (a) The reconstructed model; (b) The progression of the regularization parameter, $\alpha^{(k)}$, with iteration $k$; (c) The progression of the relative error $RE^{(k)}$ at each iteration; (d) The UPRE functional at iteration $k=5$.}

This small model permits comparison of the spectrum of the  projected matrix obtained from the RSVD algorithm  with that of the FSVD. We present this comparison for a sample iteration, iteration $5$, for two values of $q$, i.e. $q=100$ and $q=200$, in Figures~\ref{fig:figure3a} and~\ref{fig:figure3b}. The singular values of the projected matrix are consistent with the $q$ large singular values of the original matrix, with only a slight deterioration for the smaller singular values. This verifies that  for $q\ll m$, $\Gtildetilde_q$ inherits the dominant spectrum and not the conditioning of $\Gtildetilde$. This is in contrast to the LSQR algorithm in which the singular values of the projected matrix tend to approximate both large and small singular values of the original matrix, hence inheriting the conditioning of the original matrix. This is further illustrated in    Figures~\ref{fig:figure3a} and~\ref{fig:figure3b} for the singular values of the LSQR  projected matrix for $t=100$ and $t=200$. Vatankhah et al. \shortcite{VRA:2017} showed that these small singular values impact the estimation of the regularization parameter, which is then underestimated. Truncating the spectrum and using the TUPRE mitigates the issue and a reasonable $\alpha$ is found. 

\multiplot{2}{figure3a,figure3b}{height=0.20\textheight}
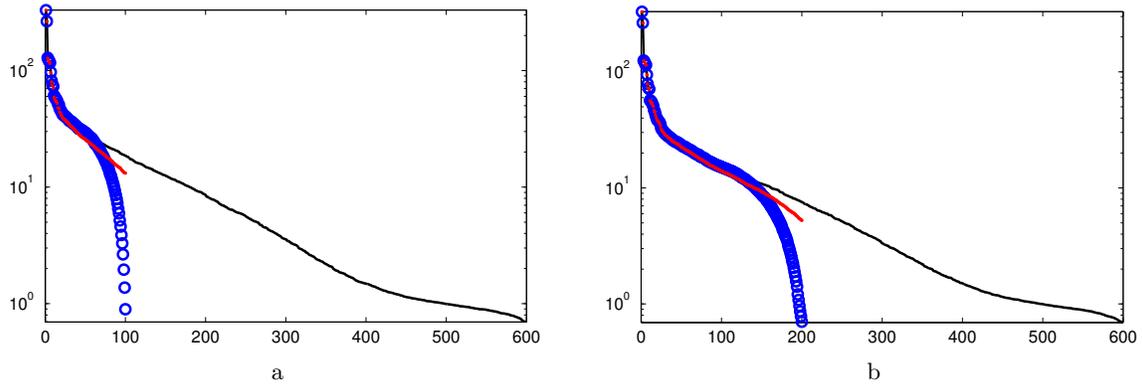
{Comparing the singular values of the original matrix, indicated by black line, and the singular values of the projected matrix, indicated by red $\cdot$ for the RSVD algorithm and blue $\circ$ for the LSQR algorithm, at iteration $5$. The projected matrix is obtained in (a) For RSVD with $q=100$ and LSQR with $t=100$ and (b) For RSVD with $q=200$ and LSQR with $t=200$.}

\subsection{Model of multiple bodies}

We now apply the inversion methodology on a larger model consisting of six bodies with different shapes and dimensions. Figure~\ref{fig:figure4a} shows a perspective view of the model, and  four plane-sections of the model are illustrated in Figure~\ref{fig:figure5a,figure5b,figure5c,figure5d}. The gravity data of the model was generated at $100 \times 55 =5500$ points on the surface with spacing $50$~m.  The noise with standard deviation of $ (0.02~(\bfde)_i + 0.001~\|\bfde\|)$ was added to provide noisy data $\bfdo$, Figure~\ref{fig:figure4b}. For the inversion, the subsurface is divided into $ 100 \times 55 \times 12 =66000$ cells of size  $50$~m in each dimension. For the inversion methodology presented in Algorithm~\ref{IterativeL1RSVD} we use $\bfma=\mathbf{0}$, $\rho_{\mathrm{min}}=0$~g~cm$^{-3}$, $\rho_{\mathrm{max}}=1$~g~cm$^{-3}$ and $K_{\mathrm{max}}=50$. We present the results of the inversion for both RSVD and LSQR algorithms for different $q$ and $t$, respectively, in Table~\ref{tab2}, noting that the problem is large and the FSVD is not feasible computationally. Except for very small $q$ it is clear that the  RSVD methodology yields acceptable solutions with relative errors that are close to those obtained using the LSQR but at much reduced CPU time.  We note that in Algorithm~\ref{RSVDAlgorithm}, the original matrix $\Gtildetilde$ is visited only twice, while the LSQR requires $t$  multiplications with $\Gtildetilde$ and $(\Gtildetilde)^T$. Furthermore, it is necessary to use reorthogonalization when using the LSQR. The inversion methodology based on the LSQR is far more expensive than that using the RSVD, which is especially evident as the problem size increases.  To illustrate the inversion results, the plane-sections of the reconstructed model using Algorithm~\ref{IterativeL1RSVD} for case $q=1000$ are shown in Figure~\ref{fig:figure6a,figure6b,figure6c,figure6d}. The recovered model is in good agreement with the original model at shallow to intermediate depths. The extent of the some of the bodies is overestimated but the horizontal borders are recovered well. An isosurface of the solution, the regularization parameter and the relative error at each iteration, and UPRE functional at the final iteration are presented in  Figures~\ref{fig:figure7a,figure7b,figure7c,figure7d}. 

Now the parameter $q$ determines the dimension of the projected subspace and it is essential that $q$ is chosen to control both the accuracy and the efficiency of the Algorithm~\ref{RSVDAlgorithm}.  Most crucially, $q$ should be large enough that the dominant spectrum of $\Gtildetilde$ is captured so that the RSVD solution carefully approximates the FSVD solution. Simultaneously, $q$ should not be so large that the computational cost becomes prohibitive. Our investigation of the gravity inverse problem, which as noted is only mildly ill-conditioned,  suggests that  $q \geq (m/6)$ provides a good compromise in using Algorithm~\ref{RSVDAlgorithm}. For problems which are severely ill-conditioned it would be feasible to use a smaller $q$.

\multiplot{2}{figure4a,figure4b}{width=0.4\textwidth}
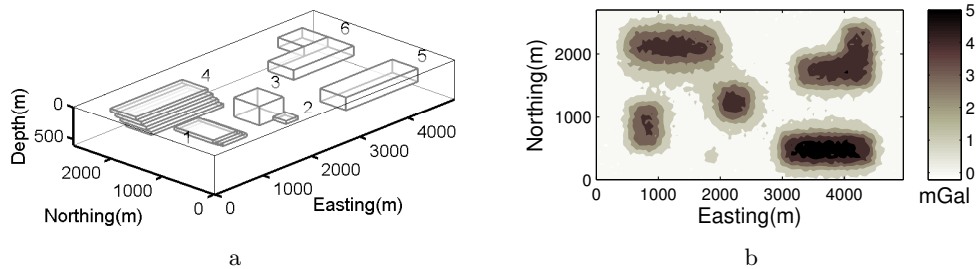
{(a) Perspective view of a model consisting of six bodies with different shapes and dimensions, bodies have density contrast $1$~g~cm$^{-3}$. (b) The data of the model generated at the surface and contaminated with noise.} 

\multiplot{4}{figure5a,figure5b,figure5c,figure5d}{width=0.4\textwidth}
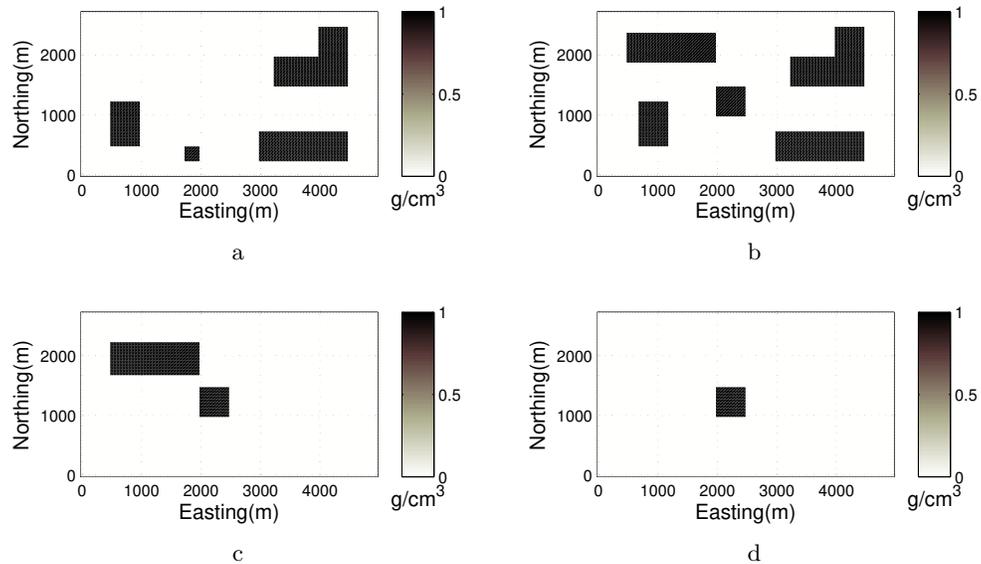
{The plane-sections of the model presented in Figure~\ref{fig:figure4a} at depths (a) $50$~m; (b) $150$~m; (c) $250$~m; (d) $350$~m.}

\begin{table}
\begin{center}
\caption{The inversion results for a model of multiple bodies obtained from inversion methodologies: LSQR and RSVD, using UPRE for RSVD but TUPRE for LSQR. Different values of parameters $q$ and $t$ are used for LSQR and RSVD, respectively.}
\label{tab2}
\begin{tabular}{c  c  c  c  c c}
\hline
Method&  $RE^{(K)}$&   $\alpha^{(1)}$ & $\alpha^{(K)}$& $K$ & Time (s)  \\ \hline
LSQR  &  & &  & &\\
$t=300$ & 0.6908& 38590& 10.43& 10& 1145.5\\ 
$t=600$ & 0.6641& 53470& 9.23& 10& 2721.6\\ 
$t=900$ & 0.6550& 64534& 10.45&  10& 4091.4\\
$t=1000$ & 0.6526& 67723& 10.62& 10 & 4642.9 \\  \hline
RSVD  &  & &  & &\\
$q=300$ & 0.7641& 46684& 2.92& 50& 334.1\\ 
$q=600$ & 0.6896& 62172& 6.05& 36& 435.2\\ 
$q=900$ & 0.6566& 72659& 13.47& 11& 202.3\\ 
$q=1000$ & 0.6556& 75604& 11.39& 10& 219.4\\  \hline
\end{tabular}
\end{center}
\end{table}

\multiplot{4}{figure6a,figure6b,figure6c,figure6d}{width=0.4\textwidth}
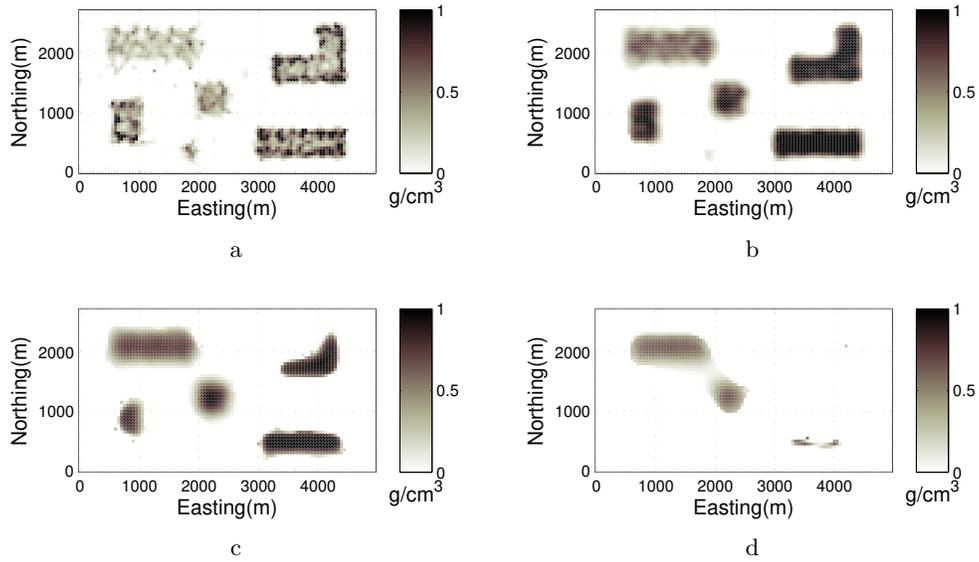
{The plane-sections of the reconstructed model for data in Figure~\ref{fig:figure4b} using Algorithm~\ref{IterativeL1RSVD} with $q=1000$. The sections are at depths (a) $50$~m; (b) $150$~m; (c) $250$~m; (d) $350$~m.}

\multiplot{4}{figure7a,figure7b,figure7c,figure7d}{width=0.4\textwidth}
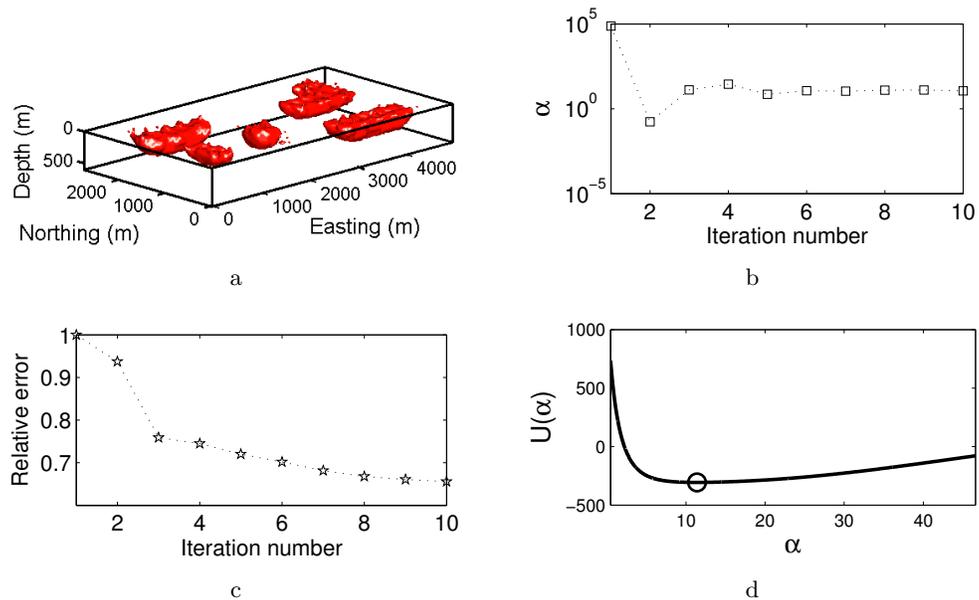
{(a) Isosurface of the reconstructed model with the density greater than $0.5$~g~cm$^{-3}$; (b) The progression of the regularization parameter, $\alpha^{(k)}$, with iteration $k$; (c) The progression of the relative error $RE^{(k)}$ at each iteration; (d) The UPRE functional at the final iteration.}

\section{Real data}
To illustrate the relevance of the approach for real data, we use gravity data over the Morro do Engenho (ME) complex in the Goi\'as Alkaline Province (GAP), in the center of Brazil. The GAP is characterized by ultramafic to felsic plutonic bodies in the north and extensive kamafugite lava flows in the central and southern part \citep{DuMa:2009}. We selected an area consisting of the two bodies in the northern part of GAP, in which one of them, Morro do Engenho (ME), is outcropping and the other, A$2$, is a possible alkaline completely buried by Quaternary sediments \citep{DuMa:2009}. We digitized the residual data carefully from Figure~$3$ in Dutra and Marangoni \shortcite{DuMa:2009}, see Figure~\ref{fig:figure8}. Some strong anomalies are related to outcrops and for others there is no  observable geological evidence  \citep{DuMa:2009}. A detailed geology of the area and the measurement data are presented in Dutra and Marangoni \shortcite{DuMa:2009} and we refer the readers to this paper and references therein. We select this data set because there are inversion results  using algorithms presented in Li and Oldenburg \shortcite{LiOl:98} and in Silva and Barbosa \shortcite{SiBa:2006}  which thus permit comparison with our inversion results, see Dutra and Marangoni \shortcite{DuMa:2009}.

The data was digitized into a grid with $45 \times 53= 2385$ data points with spacing $1$~km. We suppose each datum has an error with standard deviation $ (0.03~(\bfdo)_i + 0.003~\| \bfdo \|)$. For the inversion we use a model consisting of cells with dimension of $1$~km. For the model extending to depth $14$~km,  there are $45 \times 53 \times 14= 33390$ model parameters. Based on geological information, following Dutra and Marangoni \shortcite{DuMa:2009}, density limits $\rho_{\mathrm{min}}=0$~g~cm$^{-3}$ and $\rho_{\mathrm{max}}=0.3$~g~cm$^{-3}$ are imposed. We start the inversion with an initial model in which for nine cells in the  first layer the density contrast is selected as $0.3$~g~cm$^{-3}$  and other model parameters are zero. These cells are located at the outcrop, and, as mentioned in the discussion on the inversion methodology, the corresponding entries on the diagonal of $\Wh$  are selected to be large, here $100$, to force the algorithm to maintain the selected initial density for these cells during the iterations. We use Algorithm~\ref{IterativeL1RSVD} with $q=400$. 

The inversion terminated after $11$ iterations, using just  $34$ seconds of CPU time. The reconstructed model is shown in six plane-sections in Figure~\ref{fig:figure9a,figure9b,figure9c,figure9d,figure9e,figure9f}. These sections are consistent with those illustrated by Dutra and Marangoni \shortcite{DuMa:2009} in figure~$7$. In the first layer, the density for the selected nine cells over the outcrop was kept fixed. There is no other significant anomaly in this layer, which is consistent with the known geology of the area. The anomaly ME extends from the surface to a depth of $10-11$~km, while A$2$ starts from $2-3$~km and extends to $8$~km. There is a connection between the two anomalies at depths from $4$~km to $6-7$~km. The results are close to those obtained in Dutra and Marangoni \shortcite{DuMa:2009} and indicate the effectiveness of the RSVD algorithm. Differences are evident at the borders of the reconstructed anomalies. This is a feature of the focusing algorithm used here, as compared to the results in Dutra and Marangoni \shortcite{DuMa:2009} that are smooth. The main advantage of the RSVD algorithm is its computational time. Larger problems will be solved with acceptable computational costs. We illustrate the progression of the regularization parameter at each iteration in Figure~\ref{fig:figure10a} and the UPRE functional at the final iteration in Figure~\ref{fig:figure10b}.

\plot{figure8}{height=0.15\textheight}
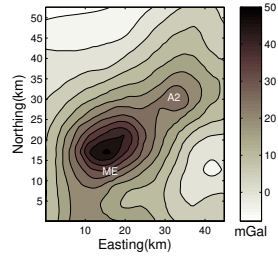
{Residual gravity data over the Morro do Engenho complex of Central Brazil. The data were digitized from figure~$3$ in Dutra and Marangoni \shortcite{DuMa:2009}.}

\multiplot{6}{figure9a,figure9b,figure9c,figure9d,figure9e,figure9f}{height=0.15\textheight}
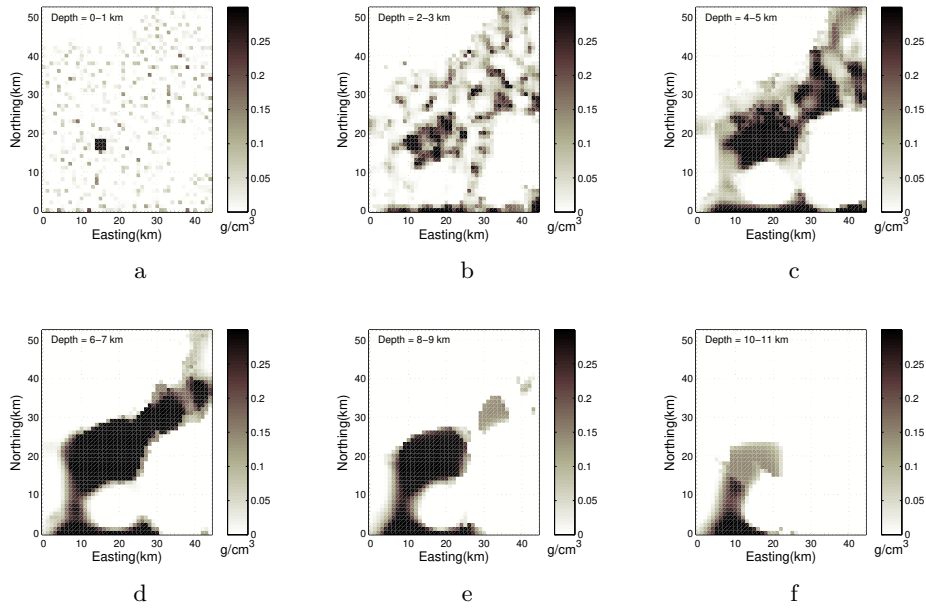
{The plane-sections of the reconstructed model for the data in Figure~\ref{fig:figure8} using Algorithm~\ref{IterativeL1RSVD} with $q=400$. The sections are at the depths specified in the figures. }   

\multiplot{2}{figure10a,figure10b}{height=0.15\textheight}
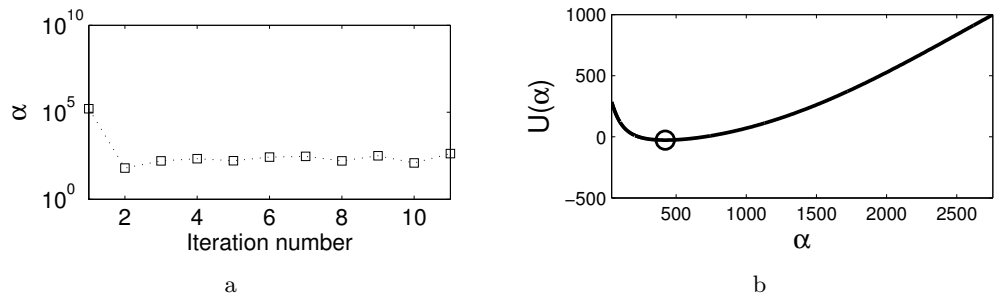
{(a) The progression of the regularization parameter, $\alpha^{(k)}$, with iteration $k$; (b) The UPRE functional at the final iteration.}

\section{Conclusion}

We have presented a fast algorithm for $3$-D inversion of gravity data based on the use of the randomized singular value decomposition. At the heart of the presented inversion methodology is a new algorithm for low-rank SVD approximation of a large under-determined matrix that makes it feasible to compute the large singular values of the original matrix in a short time. The dominant computational cost of the new RSVD algorithm is $O(lmn)$.  Our analysis showed that if a low rank $q$ approximation of the original matrix  is selected such that most of the dominant singular values are approximated, then the error between the solutions obtained using the  RSVD and the full SVD is very small. Furthermore, we demonstrated that the UPRE parameter choice rule can be used for the projected space without requiring truncation of small singular values. This is a significant difference, besides the CPU time, with the inversion methodology based on the LSQR algorithm. We showed the efficiency of the presented inversion methodology using different synthetic tests and a real case data from the Morro do Engenho (ME) complex in the Goi\'as Alkaline Province (GAP) in the center of Brazil.

\section{ACKNOWLEDGMENTS}
Rosemary A. Renaut  acknowledges the support of NSF grant DMS 1418377: ``Novel Regularization for
Joint Inversion of Nonlinear Problems".

\append[eigtosv]{Obtaining singular values and vectors from eigenvalues and vectors}

The reduced SVD of matrix $\B \in \mathcal{R}^{m \times l}$, step $4$ of the Algorithm~\ref{RSVDAlgorithm}, is given by $\B= \U_l \Si_l \tilde{\V}_l^T$, \citep{TrBa:97}. Here, $\U_l \in \mathcal{R}^{m \times l} $ and $ \tilde{\V}_l \in \mathcal{R}^{l \times l}$ are the left and the right singular vectors of $\B$ with columns denoted by $\bfui$ and $\tilde{\mathbf{v}}_i$, respectively. The diagonal matrix $\Si_l \in \mathcal{R}^{l \times l}$ contains singular values of the $\B$ ordered as $\sigma_1 \geq \sigma_2 \geq \dots \geq \sigma_{l} \geq0$. Furthermore, $ \Si_l$ and $\U_l$ are the singular values and left singular vectors of the low-rank approximation of $\Gtildetilde$, i.e. $\Gtildetilde_l$, see Xiang \& Zou \shortcite{XiZo:2013} and Voronin et. al \shortcite{VMN:2015}. The right singular vectors of the $\Gtildetilde_l$ are obtained via $\V_l= \bfQ \tilde{\V}_l \in \mathcal{R}^{n \times l}$ \citep{XiZo:2013}. Now for $\B^T\B =\tilde{\V}_l \Si_l^T \U_l^T \U_l \Si_l \tilde{\V}_l^T$ with $\U_l^T \U_l=I$, Voronin et. al \shortcite{VMN:2015}, it is immediate that $\B^T\B = \tilde{\V}_l \D \tilde{\V}_l^T $, where $\D=\Si_l^T \Si_l$. This indicates that the eigen-decomposition of  $\B^T\B$ gives the singular values, $\Si_l= \sqrt{\D_l}$, and first $l$ right singular  vectors of the $\B$, $\tilde{\V}_l$.  Finally, to compute the left singular vectors we note that $\B \tilde{\V}_l= \U_l \Si_l \tilde{\V}_l^T \tilde{\V}_l = \U_l \Si_l$, then $\U_l=\B \tilde{\V}_l \Si_l^{-1} $. In this way, we can avoid computing the SVD of matrix $\B$ directly, instead we use the eigen-decomposition of the smaller matrix $\B^T\B$. We note also, that in the generation of the rank $q$ approximation it is only the dominant $q$ terms of the spectral decomposition that are required, and hence all terms with $l$ columns can be replaced by those with just $q$ columns.

\newpage

\bibliographystyle{seg}  

\end{document}